\documentclass{amsart}

\usepackage{graphicx}
\usepackage{eepic}
\usepackage{epic}
\usepackage{amssymb,latexsym}
\usepackage{amsmath}
\usepackage{amsthm}
\usepackage{amssymb}
\usepackage{color, soul}
\usepackage{hyperref}
\usepackage{float}
\restylefloat{figure}
\restylefloat{table}
\newcommand{\cOSC}{OSC${}_{\mathrm{co}}$}
\newcommand{\cFTC}{FTC${}_{\mathrm{co}}$}
\renewcommand{\int}{\mathrm{int}}

\newtheorem{theorem}{Theorem}[section]

\newtheorem{corollary}[theorem]{Corollary}

\DeclareMathOperator{\hull}{hull}

\title[WSP $\neq$ \cFTC]{Equicontractive weak separation property on the 
    line does not imply convex finite type condition}
\author{Kevin G. Hare}
\address{Department of Pure Mathematics \\
          University of Waterloo \\
          Waterloo, Canada}
\email{kghare@uwaterloo.ca}
\thanks{Research of K.G. Hare was supported, in part, by NSERC Grant 2019-03930}

\subjclass[2020]{28A80}
\keywords{Iterated Function System, Weak Separation Property, Finite Type Condition}

\begin{document}
\begin{abstract}
Let $\{S_1, S_2, \dots, S_n\}$ be an iterated function system on $\mathbb{R}$ 
    with attractor $K$.
It is known that if the iterated function system satisfies the weak separation
    property and $K = [0,1]$ then the iterated function system also 
    satisfies the convex finite type condition.
We show that the condition $K = [0,1]$ is necessary. 
That is, we give two examples of 
    iterated function systems on $\mathbb{R}$ satisfying
    weak separation condition, and $0< \dim_H(K) < 1$ such that the IFS
    does not satisfy the convex finite type condition.
\end{abstract}
\maketitle


\section{Introduction}
\label{sec:intro}

Let $\mathcal{S} = \{S_1, S_2, \dots, S_n\}$ be a set of linear contractions.
It is well know \cite{Edgar,Falconer, Hata} 
    that there exists a unique non-empty compact set 
    $K$ such that $K = \cup S_i(K)$.
In this case we call $\{S_1, S_2, \dots, S_n\}$ an {\em iterated function 
    system (IFS)} and we call $K$ the {\em attractor} of the IFS.

In this paper we will assume $S_i: \mathbb{R} \to \mathbb{R}$.
We further assume the maps are equicontractive in the positive direction.
That is, there exists an $0 < r < 1$ such that
    for all $i$, $S_i(x) = r x + d_i$.

One common question in the study of IFSs is, under what conditions can we 
    exactly compute the Hausdorff dimension of $K$.

The first such condition is called the {\em strong separation property (SSP)}.
We say $\mathcal{S}$ satisfies the SSP if 
    \[ \min_{i \neq j} \left(\inf_{\substack{x \in S_i(K),\\ y \in S_j(K)}} |x -y|\right) > 0. \]
Under this condition all images of $K$ are very well separated.
The dimension can be computed exactly as the unique $s > 0$ 
    such that $\sum_{i=1}^n r_i^s = 1$.
In our restricted setting, where all contraction ratios are the same, this 
    gives the dimension as $-\frac{\log(n)}{\log(r)}$.

The next such properties are the {\em open set condition, (OSC)} and the 
    {\em convex open set condition (\cOSC)}.
We say $\mathcal{S}$ satisfies the OSC if there exists an non-empty open 
    set $V$ such that $S_i(V) \subset V$ for all $i$, and $S_i(V) \cap S_j(V) = \emptyset$ for all $i \neq j$.
We say $\mathcal{S}$ satisfies \cOSC\ if it satisfies OSC with 
    $V = \int(\hull(K))$.
It is easy to see that if $K$ satisfies SSP then it satisfies OSC.

Thi next conditions considered is the {\em weak separation property (WSP)}.
We will define this for the equicontractive case only.  
Care needs to be taken when adapting this to a non-equicontractive setting, 
    or higher dimensions.
Let $\sigma, \tau \in \Sigma^k = \{1,2,\dots, n\}^k$.
That is $\sigma = s_1 s_2 \dots s_k$ and $\tau = t_1 t_2 \dots t_k$ where
    $s_i, t_i \in \{1, 2, \dots, n\}$.
We write $S_{\sigma} = S_{s_1} \circ S_{s_2} \circ \dots \circ S_{s_k}$, 
    and similarly for $S_\tau$.
We see that $S_{\sigma}^{-1} \circ S_{\tau}(x) = x + a_{\sigma, \tau}$ for
    some real number $a_{\sigma, \tau}$.
We say $\mathcal{S}$ satisfies the WSP if 
    \[ 
    \inf_{\substack{\sigma, \tau \in \Sigma^k \\ 
                    k \geq 0 \\ 
                    a_{\sigma, \tau} \neq 0}} |a_{\sigma, \tau}| > 0. \]
It is clear that if $\mathcal{S}$ satisfies \cOSC, then it satisfies 
    WSP.
In this case the above infimum would be bounded below by $|K|$, the diameter 
    of $K$.
A stronger statement is true.
It is known that $\mathcal{S}$ satisfies OSC if and only if it satisfies 
    WSP and there are no exact overlaps.
That is, there do not exist $\sigma \neq \tau \in \Sigma^*$ such that 
    $S_\sigma = S_\tau$.
(Here $\Sigma^* = \cup \Sigma^k$.)
See for example \cite[Theorem 4.2.11]{BSS}.

We can think of SSP as saying images of $K$ do not overlap, 
    OSC as saying images of $K$ do not overlap in a meaningful way, 
    and WSP as saying images of $K$ either overlap exactly, or overlap in 
    a way that isn't too close.
The last condition, finite type condition can be thought of as saying that 
    there are only a finite
    number of ways images of $K$ can overlap with each other.

The last separation conditions to be discussed are 
    the {\em finite type condition (FTC)} and the 
    {\em convex finite type condition (\cFTC)}.
We will define this for the equicontractive case.  
Let $V$ be a non-empty open set such that $S_i(V) \subset V$.
We say $\sigma, \tau \in \{1, 2, \dots, n\}^k$ are neighbours 
    if $S_{\sigma}(V) \cap S_{\tau}(V) \neq \emptyset$.
We define the neighbourhood type of $S_\sigma$ as 
    \[ 
    N(S_{\sigma}) = \{S_{\sigma}^{-1} \circ S_{\tau}: S_{\tau} \text{ is a neighbour of } S_\sigma\}. \]
We say $\mathcal{S}$ satisfies the FTC if there exists a open set 
    $V$ for which there are a finite number of neighbourhood types.
We say $\mathcal{S}$ satisfies the \cFTC\ if it satisfies the FTC with 
    $V = \int(\hull(K))$.

It is possible to extend the definition of (convex) finite type condition 
    to (convex) generalized finite type condition.
This allows for contractions that are not equicontractive, or logarithmically 
    commensurate.
This level of generalization is not needed for this paper.
We refer the reader to \cite{HHR21, LN2, NW}.
The key point needed is that if $\mathcal{S}$ satisfies 
    (convex) finite type condition then it satisfies 
    (convex) generalized finite type condition.
It is shown in \cite{HHR21} that convex generalized finite type condition 
    is equivalent to finite neighbour condition, which again we will not 
    define here.

It is clear that if $\mathcal{S}$ satisfies \cFTC, then it satisfies 
    WSP.
The stronger statement is true.  
It was shown in \cite{Ngu02} that if $\mathcal{S}$ satisfies FTC, then 
    it satisfies WSP.
It was asked in \cite{LN2} if the converse is also true.
That is, if $\mathcal{S}$ satisfies WSP then is it true that it satisfy general finite type condition.

Despite the stronger conditions on FTC than WSP, it is not clear if 
    in fact this conditions in fact give rise to different IFS.
There are examples known in $\mathbb{R}^d$ with $d \geq 2$ which satisfy
    WSP but not FTC \cite{LN2, Zerner}.
These examples have the attractor $K$ within a hyperplane of $\mathbb{R}^d$, 
    and use rotations around this hyperplane to ensure that the example
    does not satisfy FTC.
When these maps are renormalized so that they are maps from the 
    hyperplane to the hyperplane, then these IFSs again satisfy FTC.
In the other direction, it is shown that if $\mathcal{S}$ satisfies the 
    WSP and $K = [0,1]$ then $\mathcal{S}$ satisfies the \cFTC.
Initially this was proved in \cite{Fe16} with the added restriction
    that the IFS was equicontractive, and all contractions were in the positive
    orientation.
This proof was later extended in \cite{HHR21} to allow for non-equicontractive maps,
    and contractions with negative orientation.
The proof in \cite{HHR21} still required $K = [0,1]$.
More precisely, 

\begin{corollary}[Corollary 4.6 of \cite{HHR21}]
Suppose the IFS $\mathcal{S}$ has self-similar set $[0,1]$.
Then then following are equivalent:
\begin{enumerate}
\item $\mathcal{S}$ satisfies the weak separation property;
\label{equiv:1}
\item $\mathcal{S}$ satisfies the finite neighbour condition;
\label{equiv:2}
\item $\mathcal{S}$ satisfies the convex generalized finite type condition;
\label{equiv:3}
\item There exists some $c > 0$ such that for any $0 < \alpha < 1$, words
    $\sigma, \tau \in \Lambda_{\alpha}$ and $z, w \in \{0,1\}$, 
    either $S_{\sigma}(z) = S_{\tau}(w)$ or 
    $|S_{\sigma}(z)-S_{\tau}(w)| > c \alpha$.
\label{equiv:4}
\end{enumerate}
\end{corollary}

It was then asked if this was true if we relax the hypothesis.

In this paper we give an two examples of IFSs $\mathcal{S}$ with 
    $[0,1] \neq K \subset [0,1]$ and where 
    \eqref{equiv:1} and \eqref{equiv:4} are true, but
    \eqref{equiv:2} and \eqref{equiv:3} are not.
The first example, given in Section \ref{sec:ex 1} has the advantage that it is simpler to 
    construct.
It has the property that it 
    satisfies WSP, OSC, FTC and does not satisfy \cOSC\ nor \cFTC.

It is worth noting that if $\{S_1, S_2, \dots, S_n\}$ satisfies the OSC, 
    then $\{S_1, S_1', S_2, \dots, S_n\}$ with $S_1 = S_1'$ 
    will satisfy FTC, but will not satisfy OSC.
This is because there is an exact overlap of level one cylinders of the maps.
This is a trivial reason, and no of interest as an example.

The second example, given in Section \ref{sec:ex 2} is a more 
    complicated construction. 
It has the property that it 
    satisfies WSP, FTC, and does not satisfy OSC, \cOSC\ nor \cFTC.
There is an exact overlap of level two cylinders, but no exact overlap 
    of level one cylinders.  

Both examples are equicontractive and satisfy $0, 1 \in K \subset [0,1]$.

\section{Construction of first example}
\label{sec:ex 1}

\subsection{Overview of IFS} {\ }

The constructed IFS will be of the form 
\[ \mathcal{S} = \{S_1, S_2, S_3\} \]
where 
\begin{align*}
S_1(x) & = x/7 \\    
S_2(x) & = x/7 + a\\    
S_3(x) & = x/7 + \frac{6}{7}\\    
\end{align*}
Here $a$ is chosen such that the IFS satisfies the WSP and does not satisfy \cFTC.
One choice gives an approximate value of $a$ as
\begin{align}
a & \approx \frac{0.9482520978}{7} \approx 0.1354645854 \label{eq:a1}
\end{align}
A precise description of $a$ will be given in Section \ref{ssec:OSC example}.

\subsection{Constructing the IFS} {\ }

As we wish to show that the IFS does not satisfy \cFTC, and 
    $0, 1 \in K \subset [0,1]$ we may assume the open set in the 
    definition of FTC is $(0,1)$.
Let $\sigma^{(1)} = 1$ and $\tau^{(1)} = 2$.
We see that $S_{\sigma^{(1)}}(0) < S_{\tau^{(1)}}(0) < S_{\sigma^{(1)}}(1)$
    if and only if $a \in J_1 := (0, 1/7)$.
We see in this case that $N(S_{\sigma^{(1)}}) = \{\mathrm{Id}, S_{\sigma^{(1)}}^{-1} \circ S_{\tau^{(1)}}\}$ 
    is a non-trivial neighbourhood type with two elements.

Assume what have a $\sigma^{(n)}$, $\tau^{(n)}$ and 
    $J_n$ such that 
    $S_{\sigma^{(n)}}(0) < S_{\tau^{(n)}}(0) < S_{\sigma^{(n)}}(1)$
    if and only if $a \in J_n$.
We will extend these to a $\sigma^{(n+1)}$, $\tau^{(n+1)}$ and $J_{n+1}$ 
    such that these properties continue to hold.
Further, we will ensure that for all $a \in J_{n+1}$ and all $\sigma \neq \tau$ of length $n$ that 
    $|S_{\sigma}^{-1}\circ S_{\tau}| \geq 4/7$.
This will show that for $a \in \cap J_n$ that the IFS satisfies the WSP.
As the $J_n$ are open sets, isn't immediately clear that $\cap J_n$ is 
    non-empty.
This will be addressed later.

We extend $\sigma^{(n)}, \tau^{(n)}$ and $J_n$ in one of two different ways 
    at each step to ensure these properties continue to hold.
We then exploit the fact that we have an infinite number of choices to ensure 
    that $\cap J_n$ is non-empty, and such that for $a \in \cap J_n$ we have that $N(S_{\sigma^{(n)}})$ are all distinct.

Assume we have such a $\sigma^{(n)}, \tau^{(n)}$ and $J_n$.
We see that both $S_{\sigma^{(n)}}(0)$ and $S_{\tau^{(n)}}(0)$ are linear functions 
    in $a$.
Further, for $\sigma^{(n)} = s_1 s_2 \dots s_n$, we see that the slope 
    of this function is $\sum_{i: s_i = 1} \frac{1}{7^i}$.
We create $\sigma^{(n)}$ and $\tau^{(n)}$ by extending $\sigma^{(n-1)}$ and $\tau^{(n-1)}$, which
    in turn were initially extended from $\sigma^{(1)} = 1$ and $\tau^{(1)} = 2$.
(We occasionally reverse the roles of $\sigma^{(k)}$ and $\tau^{(k)}$ so we do not
    know which one was an extension of $\sigma^{(1)}$.)
As one of $\sigma^{(n)}$ and $\tau^{(n)}$ has initial term 1 and one of them has initial term 2 we see 
    $S_{\sigma^{(n)}}(0) - S_{\tau^{(n)}}(0)$ is a non-constant linear function
    with respect to $a$.
Denote this as $T_n(a) = S_{\sigma^{(n)}}(0) - S_{\tau^{(n)}}(0)$.

We will choose $J_{n+1}$ in one of two ways.
The first option is to choose $J_{n+1}$ such that 
       $a \in J_{n+1}$ if and only if $S_{\sigma^{(n)} 3}(0) < S_{\tau^{(n)}1} (0) < S_{\sigma^{(n)} 3}(1)$.
In this case we would set $\sigma^{(n+1)} = \sigma^{(n)} 3$ and $\tau^{(n+1)} = \tau^{(n)} 1$.

The second option is to choose $J_{n+1}$ such that 
       $a \in J_{n+1}$ if and only if $S_{\tau^{(n)} 2}(0) < S_{\sigma^{(n)} 3}(0) < S_{\tau^{(n)} 2}(1)$.
In this case we would set $\sigma^{(n+1)} = \tau^{(n)} 2$ and $\tau^{(n+1)} = \sigma^{(n)} 3$.

We see on $J_n$ that $T_n(a) = S_{\tau^{(n)}}(0) - S_{\sigma^{(n)}}(0)$ is a 
    non-constant linear function whose image is $(0, 1/7^n)$.

To show the existence of an $J_{n+1}$ for the first option, 
    let $J_{n+1} \subset J_n$ such that $T_n(J_{n+1}) = (6/7^{n+1}, 1/7^n)$.
Note that $S_{\sigma^{(n+1)}}(0) = S_{\sigma^{(n)}3}(0) = S_{\sigma^{(n)}}(0) + 6/7^{n+1}$ and 
    $S_{\tau^{(n+1)}}(0) = S_{\tau^{(n)}}(0)$.
Hence, on this range $T_{n+1}(a) = S^{\tau^{(n)}1}(0) - S^{\sigma^{(n)}3}(0) =
    S_{\tau^{(n)}}(0) -  S^{\sigma^{(n)}}(0) - 6/7^{n+1}$ has
    image $(0, 1/7^{n+1})$.
This proves that $a \in J_{n+1}$ if and only if $S_{\sigma^{(n+1)}}(0) < S_{\tau^{(n+1)}}(0) < S_{\sigma^{(n+1)}}(1)$
    as required.

The proof the existence of $J_{n+1}$ for the second option is similar, by 
    letting $J_{n+1} \subset J_n$ be such that 
    $T_n(J_{n+1}) = (4/7^{n+1},5/7^{n+1})$.

We alternate between these options in a non-periodic way.

\subsection{The IFS satisfies WSP} {\ }

From the above we note that if $J_{n+1}$ is chosen from the first option, 
    we have for all $a \in J_{n+1}$ we have 
    $|S_{\sigma^{(n)}}^{-1} \circ S_{\tau^{(n)}}(0)| \geq 6/7$.
If instead $J_{n+1}$ is chosen as the second option, we have 
    $|S_{\sigma^{(n)}}^{-1} \circ S_{\tau^{(n)}}(0)| \geq 4/7$.
By noting $J_{n+1} \subset J_n \subset J_{n-1} \subset \dots \subset J_1$
    we see for all $a \in J_{n+1}$ and all $k \leq n$ we have
    $|S_{\sigma^{(k)}}^{-1} \circ S_{\tau^{(k)}}(0)| \geq 4/7$.
Lastly we see for all $\sigma' \neq \tau'$ with $|\sigma'| = |\tau'| \leq n$
    that either $|S_{\sigma'}^{-1} \circ S_{\tau'}(0)| \geq 1$ or 
    it is the same as $|S_{\sigma^{(k)}}^{-1} \circ S_{\tau^{(k)}}(0)|$
    for some $k \leq n$.
Assuming $\cap J_n$ is non-empty, taking $a \in \cap J_n$ gives 
    $\mathcal{S}$ satisfies the WSP.

As we choose between the two options in a non-periodic way, we see that we 
    choose second option infinitely often.
As such, we have that
    $\cap J_n = \cap \mathrm{clos}(J_n)$ is non-empty, and a singleton.
This is our value $a$.
    
\subsection{The IFS satisfies \cOSC\ but not OSC} {\ }

It is easy to see that this does not satisfy the convex open set condition.

To see that it satisfies the OSC, 
    we will construct an open set $V$ such that $\cup S_i(V) \subset S$ and $S_i(V) \cap S_j(V) = \emptyset$ for $i \neq j$.

Let $I_n = \cup_{|\sigma| = n} S_{\sigma}([0,1])$.
We see that $K = \cap I_n$.
Let $V_n = \cup_{|\sigma| = n} S_{\sigma((3/7, 4/7))}$.
Set $V = \cup V_n$.

By construction we see that $S_1(V), S_2(V), S_3(V) \subset V$.
It is easy to see that $S_1(V) \cap S_3(V) = S_2(V) \cap S_3(V) = \emptyset$.
So it remains to show that $S_1(V) \cap S_2(V)= emptyset$.

We note that $S_1(V_0) \cap S_2(V_0)= \emptyset$ for trivial reasons.
By noting that $V_n \subset I_n$ and $V_n \cap I_{n+1} = \emptyset$ for
    all $n$, we see that if ${n_1} \neq {n_2}$ then 
    $S_1(V_{n_1}) \cap S_2(V_{n_2}) = \emptyset$ 
Hence we need only check that $S_1(V_n) \cap S_2(V_n) = \emptyset$ 
    for all $n$.
Note $\sigma^{(1)} = 1$ and $\tau^{(1)} = 2$.
By our construction of $a$ we can see that 
   $S_{\sigma^{(1)}}(V_n) \cap S_{\tau^{(1)}}(V_n) = 
    S_{\sigma^{(2)}}(V_{n-1}) \cap S_{\tau^{(2)}}(V_{n-1})$.
Continuing in this way we get
   $S_{\sigma^{(1)}}(V_n) \cap S_{\tau^{(1)}}(V_n) = 
    S_{\sigma^{(n+1)}}(V_{0}) \cap S_{\tau^{(n+1)}}(V_{0})$.
The last intersection is empty.

This proves that $\mathcal{S}$ satisfies OSC.

\subsection{The IFS does not satisfy \cFTC} {\ }

As the IFS satisfies OSC, it satisfies FTC with the same $V$.
In this case it would have only one neighbourhood type, 
    namely $\{\mathrm{Id}\}$.

As we choose between the two options in a non-periodic way, 
    see that $N(S_{\sigma^{(n)}})$ will all be distinct neighbourhood sets.
(Recall this neighbourhood type is taken with respect to the open set $(0,1) = 
    \int(\hull(K))$.)
To see this we note that $\sigma^{(n)}$ and $\tau^{(n)}$ are completely determined 
    by $a$.
Let $\sigma^{(n)}$ be the length $n$ word determined by $a$, 
    and similarly for $\tau^{(n)}$.
We see that if $7^n (S_{\sigma^{(n)}}(0) - S_{\tau^{(n)}}(0)) =
     7^{m} (S_{\sigma^{(m)}}0) - S_{\tau^{(m)}}(0))$, then we get
     then we have  
     $7^{n+k} (S_{\sigma^{(n+k)}}(0) - S_{\tau^{(n+k)}}(0)) =
     7^{m+k} (S_{\sigma^{(m+k)}}0) - S_{\tau^{(m+k)}}(0))$ for all $k$.
As such, this would imply the choice between two options would be
     eventually periodic.

This proves that $\mathcal{S}$ does not satisfy \cFTC.

\subsection{Example} {\ }
\label{ssec:OSC example}

For the $a$ given in equation \eqref{eq:a1}, we choose the first option 
    or second 
    option depending on if the 
    $n$th term of the Thue-Morse sequence was $0$ or $1$.
(See for example \cite{AlloucheShallit}.)
That is, the first five choices are
    first option, second option, second option, first option, and second option.
Any non-periodic sequence would have worked.

\begin{figure}
\includegraphics[scale=0.7]{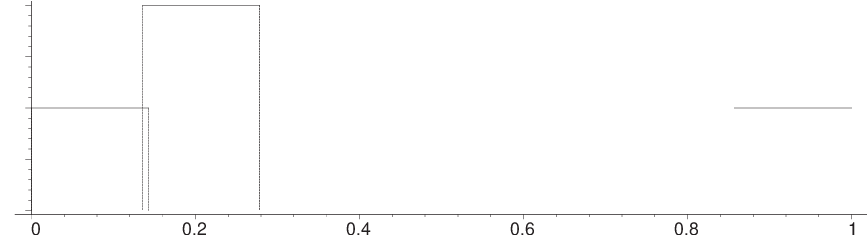} \\
\includegraphics[scale=0.7]{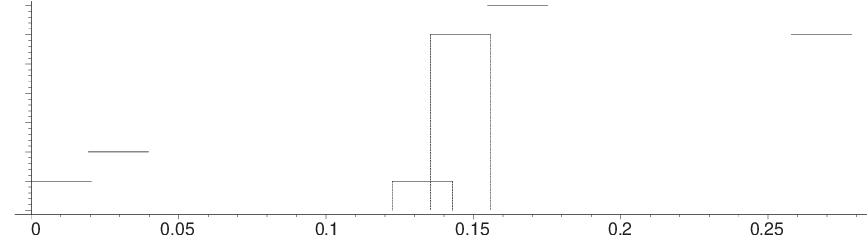} \\
\includegraphics[scale=0.7]{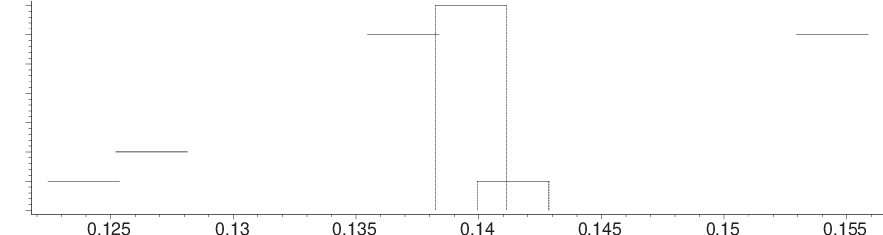} \\
\includegraphics[scale=0.7]{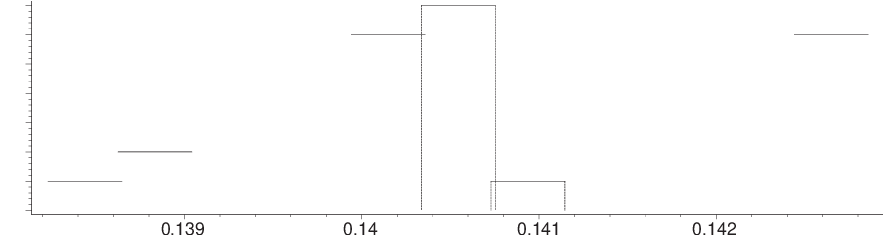} \\
\includegraphics[scale=0.7]{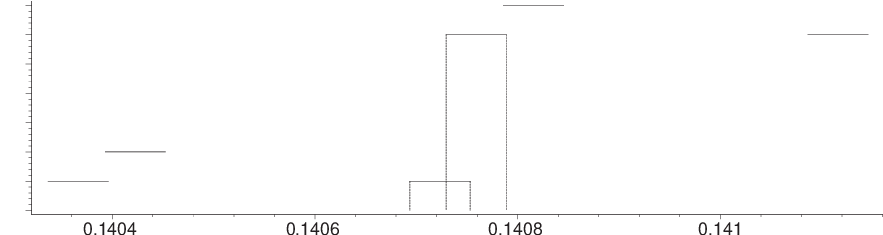} \\
\includegraphics[scale=0.7]{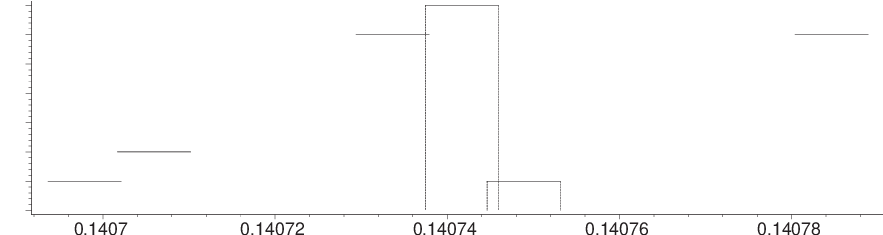}
\caption{Level $n$ cylinders for Example \ref{ssec:OSC example}}
\label{fig:OSC}
\end{figure}

In Figure \ref{fig:OSC} we present the overlap at level $n$ for each new 
    type.  

The first graph is the level 1 maps, $S_0([0,1]), S_1([0,1])$ and $S_2([0,1])$.
We note that $S_0([0,1])$ and $S_1([0,1])$ overlap, as is shown by the vertical lines.
This is the first non-trivial neighbourhood type.
Here $\sigma^{(1)} = 0$ and $\tau^{(1)} = 1$.

For the second graph, we expand to level 2 cylinders the non-trivial neighbourhood type from the first graph.
The upper three intervals are those coming from the right neighbour of the neighbourhood 
    type from the first graph, and the lower three intervals from the left most neighbour.
We see that the $S_{02}([0,1])$ overlaps $S_{10}([0,1])$, as indicated by the vertical lines.
This is the second non-trivial neighbourhood type.
It is worth noting that we also have overlaps at $S_{00}([0,1])$ with $S_{01}([0,1])$ and 
                          $S_{10}([0,1])$ with $S_{11}([0,1])$.
Neither of these are new neighbourhood types, as they are equivalent to those found 
    at the first level.
Here $\sigma^{(2)} = 02$ and $\tau^{(2)} = 10$, as we choose the first option.

We continue in this manner, expanding the new neighbourhood type found at level $n-1$ 
    to the level $n$ cylinders, and observing that there is a new neighbourhood type at level $n$.

\section{Construction of second example}
\label{sec:ex 2}

This is very similar to the first example.
We have two additional maps to force an exact overlap
    of level 2 cylinders.
This exact overlap allows us to show that it does not satisfy OSC.

The constructed IFS will be of the form 
\[ \mathcal{S} = \{S_1, S_2, S_3, S_4, S_5\} \]
where 
\begin{align*}
S_1(x) & = x/16 \\    
S_2(x) & = x/16 + a\\    
S_3(x) & = x/16 + 15/16 - 16a \\
S_4(x) & = x/16 + 11/16 \\
S_5(x) & = x/16 + 15/16 \\    
\end{align*}
Here $a$ is chosen such that the IFS satisfies the weak 
   separation property, but does not satisfy the convex finite type condition.
Further, $a$ is chosen so that the resulting IFS 
    does not satisfy OSC, as there is an exact overlap
    given by $S_{15} = S_{23}$.

One particular value of $a$ (again based on the Thue-Morse sequence) is 
    given by
\begin{align}
a & \approx \frac{0.7493705552}{16} \approx .04683565970 \label{eq:a2} \\
15/16-16a & \approx \frac{3.010071117}{16} \approx 0.1881294448 \nonumber
\end{align}

\subsection{Construction of IFS} {\ }

Let $\sigma^{(1)} = 1$ and $\tau^{(1)} = 2$.
Let $J_1 = (0, 1/16)$.
We see $a \in J_1$ if and only if 
    $S_{\sigma^{(1)}}(0) < S_{\tau^{(1)}}(0) < S_{\sigma^{(1)}}(1)$.

We set $\sigma^{(2)} = \sigma^{(1)} 4$ and 
    $\tau^{(2)} = \tau^{(1)}1$.
We choose $J_2 \subset J_1$ such that $a \in J_2$ if and only if 
    $S_{\sigma^{(2)}}(0) < S_{\tau^{(2)}}(0) < S_{\sigma^{(2)}}(1)$.

We next proceed as in Section \ref{sec:ex 1}, 

The first option is to choose $J_{n+1}$ such that 
       $a \in J_{n+1}$ if and only if $S_{\sigma^{(n)} 5}(0) < S_{\tau^{(n)}1} (0) < S_{\sigma^{(n)} 5}(1)$.
In this case we would set $\sigma^{(n+1)} = \sigma^{(n)} 5$ and $\tau^{(n+1)} = \tau^{(n)} 1$.

The second option is to choose $J_{n+1}$ such that 
       $a \in J_{n+1}$ if and only if $S_{\tau^{(n)} 2}(0) < S_{\sigma^{(n)} 5}(0) < S_{\tau^{(n)} 2}(1)$.
In this case we would set $\sigma^{(n+1)} = \tau^{(n)} 2$ and $\tau^{(n+1)} = \sigma^{(n)} 5$.

The proof is the same as before.

\subsection{The IFS satisfies WSP and not \cFTC, OSC, nor \cOSC} {\ }

As before, it is easy to see that the IFS satisfies WSP and not 
    \cFTC..
As there are exact overlaps, the IFS clearly does not satisfy OSC, nor
    \cOSC.

\subsection{The IFS satisfies FTC}

Although the IFS does not satisfy \cFTC, nor OSC, it does satisfy FTC.

To see that it satisfies the FTC, it suffices to construct a open set $V$ and
    demonstrate that there are only finitely many neighbourhood types with 
    this open set.

Let $I_n = \cup_{|\sigma| = n} S_{\sigma}([0,1])$.
We see that $K = \cap I_n$.
Let $V_n = \cup_{|\sigma| = n} S_{\sigma((7/16,8/16))}$.

Let $V = \cup V_n$.
We see by construction that $S_1(V), S_2(V), \dots, S_5(V) \subset V$.

We see that $S_{23}(V) = S_{14}(V) \subset S_2(V) \cap S_1(V)$.
Hence $N(S_{2}) = \{\mathrm{Id}, S_{2}^{-1}\circ S_{1}\}$ and 
      $N(S_{1}) = \{\mathrm{Id}, S_{1}^{-1}\circ S_{2}\}$.
There are no other overlaps at level 1, and hence 
    $N(S_3) = N(S_4) = N(S_5) = \{\mathrm{Id}\}$.

For the children under $S_1$ we have 
   $N(S_{11}) = N(S_1)$ and $N(S_{12}) = N(S_2)$.
All other children have neighbourhood type $\{\mathrm{Id}\}$.

Similarly the children under $S_2$ have neighbourhood types
    $N(S_1), N(S_2)$ or $\{\mathrm{Id}\}$.

This gives that there are three neighbourhood types, namely
\[ \{\mathrm{Id}\}, \hspace{1cm}
   \{\mathrm{Id}, S_{2}^{-1}\circ S_{1}\}, \hspace{1cm} 
   \{\mathrm{Id}, S_{1}^{-1}\circ S_{2}\}. \]
which proves that this IFS satisfies FTC.

\subsection{Example}{\ }
\label{ssec:FTC example}

\begin{figure}
\includegraphics[scale=0.7]{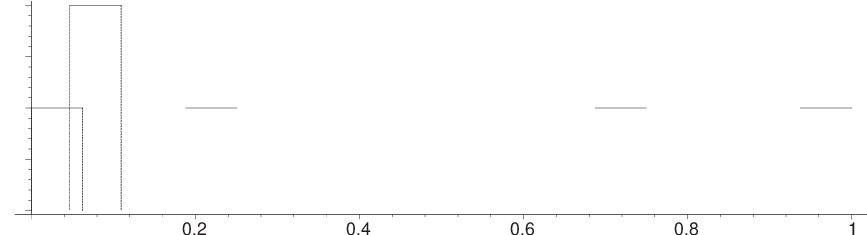} \\
\includegraphics[scale=0.7]{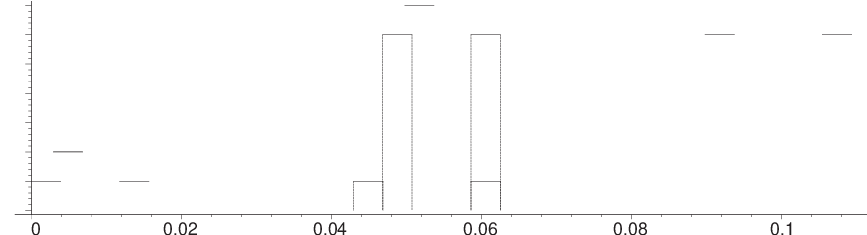} \\
\includegraphics[scale=0.7]{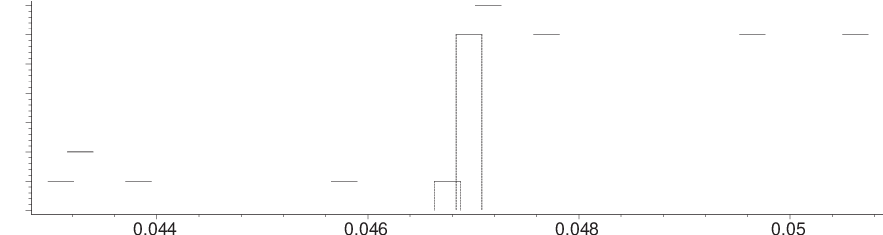} \\
\includegraphics[scale=0.7]{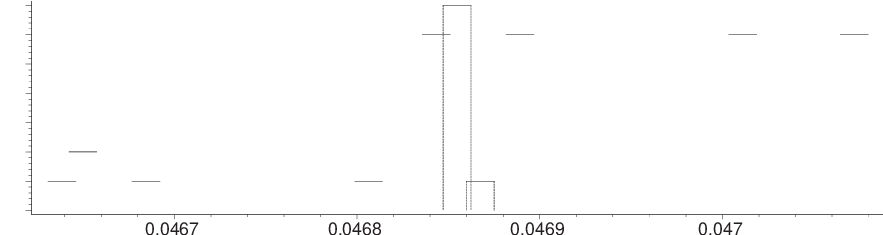} \\
\includegraphics[scale=0.7]{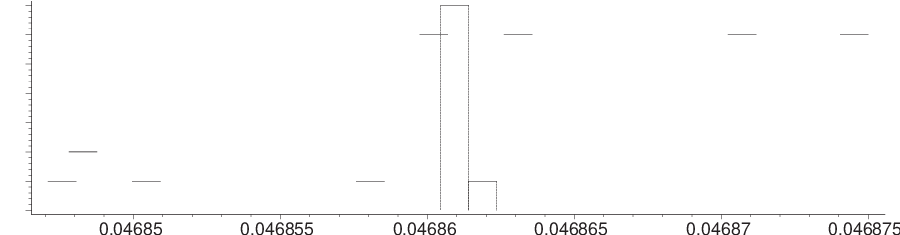} \\
\includegraphics[scale=0.7]{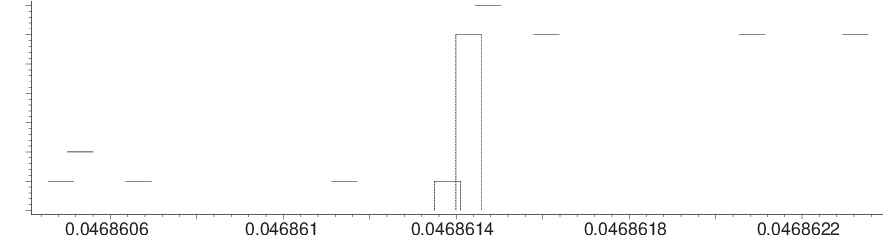}
\caption{Level $n$ cylinders for Example \ref{ssec:FTC example}}
\label{fig:NOSC}
\end{figure}

In Figure \ref{fig:NOSC} we present the overlap at level $n$ for each new 
    type.  
We choose, in order, option 1, option 2, option 2, option 1, option 2, ....
We choose this based on the Thue-Morse sequence, which is a well known
    non-periodic sequence, although any non-periodic sequence would have 
    worked.

\end{document}